\documentclass[reqno]{amsart}          
\usepackage{epsfig}                     
\usepackage{amscd}
\usepackage{amssymb}
\usepackage[all]{xy}
\usepackage{amsmath}
\usepackage{amsthm}
\newcommand{\tuborg}{\left\{\begin{array}{ll}}
\newcommand{\sluttuborg}{\end{array}\right.}

\newcommand{\calZ}{\mathcal{Z}}
\newcommand{\calD}{\mathcal{D}}
\newcommand{\calM}{\mathcal{M}}

\newcommand{\calP}{\mathcal{P}}

\newcommand{\bbZ}{\mathbb{Z}}
\newcommand{\bbP}{\mathbb{P}}
\newcommand{\bbQ}{\mathbb{Q}}
\newcommand{\bbR}{\mathbb{R}}
\newcommand{\bbC}{\mathbb{C}}

\newcommand{\bbA}{\mathbb{A}}

\newcommand{\ord}{{\rm ord}}

\newcommand{\supp}{{\rm Supp}}

\newcommand{\spec}{{\rm Spec}}

\newcommand{\res}{{\rm res}}

\newcommand{\Li}{{\rm Li}}

\newtheorem{thm}{Theorem}[section]

\newtheorem{lemma}[thm]{Lemma}
\newtheorem{cor}[thm]{Corollary}
\newtheorem{prop}[thm]{Proposition}

\newtheorem{remark2}[thm]{Remark}

\DeclareSymbolFont{AMSb}{U}{msb}{m}{n}
\DeclareMathSymbol{\N}{\mathbin}{AMSb}{"4E}
\DeclareMathSymbol{\Z}{\mathbin}{AMSb}{"5A}
\DeclareMathSymbol{\R}{\mathbin}{AMSb}{"52}
\DeclareMathSymbol{\Q}{\mathbin}{AMSb}{"51}
\DeclareMathSymbol{\I}{\mathbin}{AMSb}{"49}
\DeclareMathSymbol{\C}{\mathbin}{AMSb}{"43}

\begin{document}
 
\title{Algebraic cycles and additive dilogarithm}
\author{Jinhyun Park}
\address{Department of Mathematics, Purdue University, West Lafayette, Indiana USA}
\email{jinhyun@math.purdue.edu}
\date{July 14, 2007}
\begin{abstract} For an algebraically closed field $k$ of characteristic $0$, we give a cycle-theoretic description of the additive $4$-term motivic exact sequence associated to the additive dilogarithm of J.-L. Cathelineau, that is the derivative of the Bloch-Wigner function, via the cubical additive higher Chow groups under one assumption. The $4$-term functional equation of Cathelineau, an additive analogue of Abel's $5$-term functional equation, is also discussed cycle-theoretically.
\end{abstract}
\maketitle

\section*{Introduction}

For a complex number $z$ with $|z|<1$, the power series expansion of the function $\Li_1 (z):= - \log (1-z)= \sum_{k=1} ^{\infty} \frac{z^k}{k}$ encourages one to define the $n$-logarithm function as $\Li_n (z) := \sum_{k=1} ^{\infty} \frac{z^k}{k^n}.$ After analytic continuation one can see them as multi-valued meromorphic functions on $\bbC$. Various speculations on them suggest their connections to some arithmetic questions. For instance, in connection with the Beilinson regulators in \cite{Be}, it is believable that for a smooth complex variety $X$, the $n$-th Chern character maps $$ch_n : K_m (X) \to H_{\calD} ^{2n-m} (X, \bbR (n))$$ from the Quillen $K$-groups to the Deligne-Beilinson cohomology groups may be expressible in terms of these polylogarithm functions. When $X= \spec (\bbC)$, the generalized Bloch-Wigner functions $\calD_n (z)$ (see \cite{B3} for the original definition, and \cite{C} for generalizations), the single-valued real analytic cousins of $\Li_n (z)$ on $\bbC - \{ 0, 1 \}$, may induce the cohomology classes giving the Borel regulator elements, and the functional equations satisfied by $\calD_n(z)$ may correspond to the cocycle conditions. For a general discussion, see \cite{Hain} and \cite{Z}.

More generally for a field $k$, these functional equations can be formally used to give some relations among generators in certain free abelian groups, and the related complexes called the polylogarithmic motivic complexes seem to capture some of the rational motivic cohomology groups for $k$ (see \cite{G1}, for example).

The basic example is Abel's $5$-term functional equation \begin{equation}\label{abel}\calD_2 (x) - \calD_2 (y) + \calD_2 \left( \frac{y}{x} \right) - \calD_2 \left( \frac{1-y}{1-x} \right) + \calD_2 \left( \frac{1- y^{-1}}{1-x^{-1}}\right) = 0\end{equation} satisfied by the Bloch-Wigner function $\calD_2 (z)$. Using this relation formally, when $k$ is an infinite field, one can construct the following $4$-term motivic exact sequence (see \cite{DS}, \cite{Suslin})
\begin{equation}\label{motivic exact sequence}0 \to K_3 ^{\rm ind} (k) \otimes \bbQ \to \calP(k)  \otimes \bbQ \to \left(k^{\times} \wedge k^{\times} \right) \otimes \bbQ \to K_2 ^{M} (k) \otimes \bbQ \to 0,\end{equation} where
$$\calP(k) = \frac{\bbZ[k \backslash \{ 0, 1 \}]}{\left( \left< a \right> - \left< b \right> + \left< \frac{b}{a} \right> - \left< \frac{1-b}{1-a} \right> + \left< \frac{1 - b^{-1}}{1- a^{-1}}\right> \right)}.$$ 

Several interesting facts related to this sequence are known. When $k= \bbC$, the function $\calD_2 (z)$ is defined on the group $\calP(\bbC)$ due to \eqref{abel} (see \cite{B3}) and this map turns out to correspond to the volume map for the scissors congruence group of the 3-dimensional real hyperbolic space.  \cite{D3} and \cite{DS} show a construction of the basic exact sequence for this scissors congruence group. This sequence looks surprisingly similar to \eqref{motivic exact sequence} (see \cite{G2} and \cite{Hain}).

On the other hand, the identification of $K_3 ^{\rm ind} (k)$ with a higher Chow group $CH^2 (k,3)$ (see \cite{B1} for the definition, \cite{BL} \cite{Suslin} for the proofs), and $K_2 ^M (k)$ with $CH^2 (k,2)$ (see \cite{NS}, \cite{T}) relate algebraic cycles to this exact sequence. In \cite{GMS}, one finds that the group $CH^2 (k,3) \otimes \bbQ$ has a family of elements satisfying the same type of functional equation as $\eqref{abel}$. 

J.-L. Cathelineau's idea (\cite{C}, \cite{EVG}) of using the derivatives of the generalized Bloch-Wigner functions leads to the infinitesimal (or, \emph{additive} in the sense of S. Bloch and H. Esnault in \cite{BE2}) polylogarithm functions, and they satisfy different functional equations. For the dilogarithm case, it takes the form
\begin{equation}\label{cathelineau identity}\left<a \right> - \left< b \right> + a \left< \frac{b}{a} \right> + (1-a) \left< \frac{1-b}{1-a} \right> = 0,\end{equation} which has first appeared in \cite{C88}. (M. Kontsevich in  \cite{EVG}, \cite{Kont} noted that the same functional equation mysteriously appears in the theory of finite polylogarithms.) When $k$ is a field of characteristic $0$, \cite{C} and \cite{G3} used this $4$-term relation formally to define the additive polylogarithmic motivic complex, and as an immediate corollary of Proposition 6 in \cite{C}, one has the additive analogue of the motivic $4$-term exact sequence of abelian groups
\begin{equation}\label{cathelineau sequence}0 \to k \to T\calP (k) \to k \otimes_{\bbZ} k^{\times} \to \Omega_{k/ \bbZ} ^1\to 0.\end{equation}
The group $T\calP(k)$ is the $k^{\times}$-module $k \left< 1 \right> \oplus \beta (k)$, where $k^{\times}$ acts on 
$$\beta (k) := \frac{k [ k \backslash \{ 0, 1 \} ]}{\left(\left< a \right> - \left< b \right> + a \left< \frac{b}{a} \right> + (1-a) \left< \frac{1-b}{1-a} \right>\right) }$$trivially, and on $k\left< n \right>:=k$ via the rule $a * v = a^{2n+1} v$ for $a \in k^{\times}$ and $v \in k$.

The main results in this paper grew out of attempts to express this $4$-term exact sequence \eqref{cathelineau sequence} and the Cathelineau identity \eqref{cathelineau identity} in terms of algebraic cycles. A $K$-theoretic version of this sequence was considered in \cite{BE2} through the localization sequence for the relative pair $(k[t], (t^2))$. Here we use the cubical additive higher Chow complex of \cite{P1}. The group of $0$-cycles $ACH_0 (k, 1)$ gives the group $\Omega_{k/\bbZ} ^1$ of absolute K\"ahler differentials (see \S 1 for the definition, \cite{BE2} for the proof), whereas the group of $1$-cycles $ACH_1 (k, 2)$ is believed to give the first group $k$ (see \cite{P1} for more details). Theorem \ref{theorem 1} proves that $ACH_1 (k, 2)$ is nontrivial using the regulator map on additive Chow groups in \cite{P1}, but its identification with $k$ seems to require further very nontrivial work. Under the assumptions that $ACH_1 (k, 2) \simeq k$ and $k$ is an algebraically closed field of characteristic $0$, Theorem \ref{dilog} will show how to construct the additive $4$-term motivic exact sequence

$$\xymatrix{0 \ar[r] & ACH_1 (k,2) \ar[r]  \ar[d]^{\simeq ?}  & T\calP ^{cy} (k) \ar[r] & k \otimes_{\bbZ} k^{\times} \ar[r] & ACH_0 (k,1) \ar[r]  \ar[d]^{\simeq} & 0 \\ & k & & & \Omega_{k/\bbZ} ^1 & }$$ together with a family of classes of $1$-cycles $C_a$ satisfying the $4$-term Cathelineau identity 
$$C_a - C_b + a* C_{\frac{b}{a}} + (1-a) * C_{\frac{1-b}{1-a}} \equiv 0$$
in the group $T\calP^{cy} (k)$ purely from algebraic cycles (see \eqref{star action} for $*$, \eqref{cycle sequence} for $T\calP^{cy}(k)$, \eqref{C_a} for $C_a$).

\section{The Additive Chow group $ACH_1 (k,2)$}

We recall some basic definitions and results on additive Chow groups from \cite{P1} without proofs. Then we will prove that the group $ACH_1 (k,2)$ is nontrivial when $k$ is a field of characteristic zero.

 \subsection{Additive Chow group and regulator}
 
Let $k$ be a field. We are concerned only about the case $X= \spec (k)$ and $m=2$ in the notations of \cite{P1} so that we will drop them from the notation. Let $V$ be a normal variety over a field $k$ and $PDiv(V)$ be the set of all prime Weil divisors on $V$. For a Weil divisor $D$ on $V$, the \emph{support} of $D$, denoted by $\supp (D)$ is the set of all prime Weil divisors $Y$ such that $\ord_Y D \not = 0$. For each $D$, the set $\supp (D)$ is finite. For Weil divisors $Y_1, \cdots, Y_n$ on $V$, the \emph{supremum} of $Y_1, \cdots, Y_n$ is a Weil divisor on $V$ defined as
$$
\sup_{1 \leq i \leq n} Y_i := \sum_{Y \in {\rm PDiv} (V)} \left( \max_{1 \leq i \leq n } \ord_Y (Y_i) \right) [Y].$$
This expression makes sense because for only finitely many $Y \in {\rm PDiv} (V)$, the number $\ord_Y (Y_i)$ is nonzero for some $1 \leq i \leq n$.

Recall the definition of the cubical additive higher Chow complex with modulus $2$ for $X= \spec (k)$. Let $$ \tuborg A := ( \bbA^1 , 2 \{ 0 \} ),\\
B := ( \square, \{ 0, \infty \} ) = (\bbP^1  - \{ 1 \} , \{ 0, \infty \} ),\\
 \Diamond_n := A \times B^n \ni (x, t_1, \cdots, t_n), \\ \widehat{\Diamond}_n:=A \times \left( \bbP^1 \right)^n. \sluttuborg $$
 
For each $i \in \{ 1, \cdots, n \}$ and $j \in \{ 0, \infty \}$, we have the codimension $1$ face maps $$\mu_i ^j : \Diamond_{n-1} \hookrightarrow  \Diamond_n  $$ $$(y,  t_1, \cdots, t_{n-1})\mapsto  (y,  t_1, \cdots, t_{i-1}, j, t_{i}, \cdots, t_{n-1})$$ and various higher codimensional face maps as well. Let $F_n \subset \Diamond_n$ be the union of the codimension $1$ faces $F_i ^j := \mu_i ^j (\Diamond_{n-1})$ for $i \in \{ 1 , \cdots, n \}$ and $j \in \{ 0, \infty \}$. For $0$-cycles, define
\begin{equation}c_0 (\Diamond_n) := \bigoplus_{\xi:\mbox{closed pt}} \bbZ \xi, \ \ \xi \in  \Diamond_n - F_n- \{x=0\}.\end{equation}

 For $l$($>0$)-dimensional cycles, define $c_l( \Diamond_n)$ inductively as follows. Suppose that $c_{l-1}(\Diamond_{n-1})$ is defined. Let 
 \begin{equation} c_l (\Diamond_n) := \underset{W}{\bigoplus} ~ \bbZ W,\end{equation} where the sum is over all $l$-dimensional irreducible closed subvariety $W \subset   \Diamond_n$ with a normalization $\nu: \overline{W} \to \widehat{\Diamond}_n$ of the Zariski-closure $\widehat{W}$ of $W$ in $ \widehat{\Diamond}_n$ satisfying the following properties:
 \begin{enumerate}
 \item $W$ intersects all lower dimensional faces properly, \emph{i.e.} in right codimensions.
 \item The associated $(l-1)$-cycle of the scheme $W \cap F_i ^j$ lies in the group $c_{l-1} ( \Diamond_{n-1})$ for all $i \in \{1, \cdots, n \}$ and $j \in \{ 0, \infty \}$. This cycle will be denoted by $\partial_i ^j W$.
 \item The following Weil divisor on the normal variety $\overline{W}$ satisfies
 \begin{equation}\label{modulus condition}
 \sup_{1 \leq i \leq n} \nu^* \{ t_i = 1 \} - 2 \nu^* \{ x = 0 \} \geq 0.
 \end{equation}
 \end{enumerate}
 Via the face maps $\partial_i ^j=\left( \mu_i ^j \right)^*: c_l (\Diamond_n) \to c_{l-1} (\Diamond_{n-1})$, we obtain the boundary map \begin{equation}\partial:= \sum_{i=1} ^n (-1)^i \left( \partial_i ^0 - \partial_i ^{\infty} \right) : c_l (\Diamond_n) \to c_{l-1} ( \Diamond_{n-1}).\end{equation} We immediately see that $\partial^2 = 0$.

 \begin{remark2}\label{premodulus2}
{\rm Note that the condition \eqref{modulus condition} is equivalent to the following: for each prime Weil divisor $Y \in \supp \left( \nu ^* \{ x = 0 \} \right)$ on $\overline{W}$, there is an index $i \in \{ 1, \cdots, n\}$ such that
\begin{equation}\label{rephrasal}\ord_Y \left( \nu^* \{ t_i = 1 \} - 2 \nu^* \{ x=0 \} \right) \geq 0.\end{equation} We write $(W, Y) \in \calM ^2 (t_i)$ if this is the case.}
\end{remark2}

Let $d_l (\Diamond_n)$ be the subgroup of $ c_l (\Diamond_n)$ generated by degenerate cycles on $ \Diamond_n$ obtained by pulling back admissible cycles on $ \Diamond_{n-1}$ via various projection maps
$$ (y,  t_1, \cdots, t_n) \mapsto (y,  t_1, \cdots, \widehat{t_i}, \cdots , t_n).$$ Define $$\calZ_l ( \Diamond_n) := \frac{c_l ( \Diamond_n)}{d_l ( \Diamond_n)}.$$We easily see that the boundary map $\partial$ on $c_{*} ( \Diamond_*)$ descends onto $\calZ_{*} ( \Diamond_*)$. This gives the cubical additive higher Chow complex for $\spec (k)$:
$$\cdots \to \calZ_3 (\Diamond _{n+1}) \overset{\partial}{\to} \calZ_2 (\Diamond_{n}) \overset{\partial}{\to} \calZ_1 (\Diamond _{n-1}) \overset{\partial}{\to} \calZ_0 (\Diamond_{n-2})\to 0.$$ Each group has a natural $k^{\times}$-action determined by the actions on $k$-rational points
\begin{equation}\label{star action}
\alpha * \left( x, t_1, \cdots, t_n \right):= \left( \frac{x}{\alpha}, t_1, \cdots, t_n \right), \ \ \ \ \ \alpha \in k^{\times},
\end{equation}and the boundary map $\partial$ is $*$-equivariant. The homology at $\calZ_l (\Diamond_{n})$ is the \emph{additive higher Chow group} $ACH_l (k, n)$. In this paper, we are primarily interested in the following piece of the cubical additive higher Chow complex
$$\cdots \to \calZ_2 (\Diamond_3) \overset{\partial}{\to} \calZ_1 (\Diamond_2) \overset{\partial}{\to} \calZ_0 (\Diamond_1) \to 0,$$ and in particular $1$-cycles play the most important roles. By our definition, for an irreducible curve $C \in \calZ_1 (\Diamond_2)$ and the normalization of its projective closure $\nu: \overline{C} \to \widehat{C}$, each prime Weil divisor $p \in \supp \left(\nu^* \{ x = 0 \} \right)$ is a closed point of $\overline{C}$ and each such $p$ satisfies $(C, p) \in \calM ^2 (t_i)$ for at least one $i \in \{ 1, 2 \}$.

A restatement of the main theorem in \cite{P1} for this special case is the following:

\begin{thm} Let $k$ be a field of characteristic $0$. Then there is a nontrivial homomorphism $R_2 :\calZ_1 (\Diamond_2)  \to k$ such that the diagram
 \begin{equation}\xymatrix{ \calZ_2 (\Diamond_3) \ar[r] ^{\partial} \ar[d] & \calZ_1 (\Diamond_2)  \ar[d]^{R_2}  \\ 0 \ar[r] & k  }\end{equation}
 is commutative. This $R_2 (C)$ for an irreducible curve $C \in \calZ_1 (\Diamond_2)$ is defined as follows:
$$R_2 (C):= \sum_{p \in \supp \left( \nu^* \{ x = 0 \} \right)} R_2 (C,p),$$ where $\nu: \overline{C} \to \widehat{C}$ is a  normalization and
$$R_2 (C,p) := \tuborg \res_p \left( \nu^* \left( \frac{1-t_1}{x^{3}} \frac{dt_2}{t_2} \right) \right),  & \mbox{ if } (C, p) \in \calM^2 (t_1), \\
 - \res_p \left( \nu^* \left( \frac{ 1-t_2}{x^{3}} \frac{dt_1}{t_1} \right) \right), & \mbox{ if } (C, p) \in \calM^2 (t_2). \sluttuborg$$For general cycles, we extend it $\bbZ$-linearly. This map induces a homomorphism
 $$R_2: ACH_1 (k,2) \to k.$$
 \end{thm}
 
 The following is convenient to compute regulator values for some concrete cycles.

\begin{prop}\label{trivial vanishing}
\begin{enumerate}
\item If $\nu^* \{ x = 0 \} = 0$, i.e. $\widehat{C}$ does not intersect $\{x =  0 \}$ in $\widehat{\Diamond}_2$, then $R_2 (C)=0$.
\item If $t_1$ or $t_2 $ is constant on $C$, then $R_2(C) = 0$.
\end{enumerate}
\end{prop}

\begin{proof} (1) is obvious, because $R_2$ is evaluated only at points lying over $\{ x = 0 \}$. For (2), suppose that for example $t_1$ is constant. Then automatically for any point $p \in \nu^* \{ x = 0 \}$, we have $(C, p) \in \calM^2 (t_2)$ and $\nu^* \omega_2 = \nu^* \left( \frac{1-t_2}{x^3} \frac{dt_1}{t_1}\right) = 0$ as $dt_1 = 0$. Thus $R_2 (C) = 0$. The other case is similar.
\end{proof}

\begin{remark2}\label{volume}{\rm The $*$-action of $k^{\times}$ has an interesting property: for $\alpha \in k^{\times}$ and $C \in \calZ_1 (\Diamond_2)$, we have
$$R_2 (\alpha*C) = \alpha^3 R_2 (C).$$ Its proof is trivial. We use this observation frequently.
}
\end{remark2}

\subsection{Cycles $C_1 ^{\cdot}$ and $C_2 ^{\cdot}$.}

Let $a, a_1, a_2 \in k$ and $b, b_1, b_2 \in k^{\times}$. Let $C_1 ^{(a_1, a_2), b}$, $C_2 ^{a, (b_1, b_2)}$ be parametrized $1$-cycles in $\calZ_1 (\Diamond_2)$ defined as follows:

$$C_1 ^{(a_1, a_2), b} = \tuborg \left\{ \left( t, \frac{(1-a_1 t)(1-a_2 t)}{1-(a_1 + a_2)t} , b \right) | t \in k \right\}, & \mbox{ if } a_1 a_2 (a_1 + a_2) \not = 0, \\ \left\{ \left( t, 1- a^2 t^2 , b \right) | t \in k \right\}, & \mbox{ if } a:=a_1 = - a_2 \not = 0, \\ 0, & \mbox{ if } a_1 a_2 = 0.\sluttuborg$$

$$C_2 ^{a, (b_1, b_2)} = \tuborg \left\{ \left( \frac{1}{a} , t, \frac{b_1 t - b_1 b_2 }{t - b_1 b_2 } \right) | t \in k \right\}, & \mbox{ if } a \not = 0, \\ 0, & \mbox{ if } a = 0.\sluttuborg$$
They are taken from \S 6 in \cite{BE2}.

\begin{lemma}\label{bilinearity} \begin{enumerate}
\item $R_2 (C_1 ^{\cdot}) = R_2 (C_2 ^{\cdot}) = 0$.
\item $$\tuborg \partial C_1 ^{(a_1, a_2), b } =  - \left( \frac{1}{a_1}, b \right) - \left( \frac{1}{a_2}, b \right) + \left( \frac{1}{a_1 + a_2} , b \right), \\

\partial C_2 ^{a, (b_1, b_2)} =  \left( \frac{1}{a}, b_1 \right) + \left( \frac{1}{a}, b_2 \right) - \left( \frac{1}{a} , b_1 b_2 \right), \sluttuborg$$
where the symbol $\left( \frac{1}{a}, b \right)$ must be interpreted as $0$ if $a = 0$.
\end{enumerate}
\end{lemma}

\begin{proof}
(1) Obviously $R_2 (C_1 ^{\cdot}) = 0$ by Proposition \ref{trivial vanishing}-(2). For $C_2 ^{\cdot}$, if $a=0$, then it is trivial. When $a \not = 0$, because $\frac{1}{a} \not = 0$, by Proposition \ref{trivial vanishing}-(1) we have $R_2 (C_2 ^{a, (b_2, b_2)}) = 0$.

\noindent (2) For $C_1 ^{\cdot}$, if $a_1 a_2 (a_1 + a_2) \not = 0$, then a direct computation
$$\tuborg \partial_1 ^0 C_1 ^{(a_1, a_2), b} = \left( \frac{1}{a_1}, b \right) + \left( \frac{1}{a_2}, b \right), \\
\partial_1 ^{\infty} C_1  = \left( \frac{1}{a_1 + a_2} , b \right), \\
\partial_2 ^{0} C_1 = 0, \\
\partial _2 ^{\infty} C_1 = 0, \sluttuborg$$
gives $\partial C_1 ^{\cdot} = - \left( \frac{1}{a_1}, b \right) - \left( \frac{1}{a_2}, b \right) + \left( \frac{1}{a_1 + a_2} , b \right).$ 

If $a = a_1 =- a_2 \not = 0$, then we see that
$$
\tuborg \partial_1 ^0 C_1 ^{(a, -a), b} = \left( \frac{1}{a}, b \right) + \left( - \frac{1}{a} , b \right), \\
\partial_1 ^{\infty} C_1 = 0, \\
\partial_2 ^0 C_1 = 0, \\
\partial_2 ^{\infty} C_1 = 0, \sluttuborg
$$ so that $\partial_1 ^{(a, -a), b} = - \left( \frac{1}{a}, b \right) - \left(- \frac{1}{a}, b \right)$. 

If $a_1 a_2 = 0$, then it is trivial.

Similarly for $C_2$, when $a \not  = 0$, we have
$$\tuborg \partial_1 ^0 C_2 ^{a, (b_1, b_2)} = \left( \frac{1}{a} , 1 \right) = 0, \\
\partial_1 ^{\infty} C_2 = \left( \frac{1}{a}, b_1 \right), \\
\partial_2 ^0 C_2 = \left( \frac{1}{a} , b_2 \right),\\
\partial_2 ^{\infty} C_2 = \left( \frac{1}{a}, b_2 b_2 \right), \sluttuborg$$ that gives $\partial C_2 = \left( \frac{1}{a}, b_1 \right) + \left( \frac{1}{a} , b_2 \right) - \left( \frac{1}{a} , b_1 b_2 \right)$, and when $a=0$ it is trivial. This proves the lemma.\end{proof}

It is interesting to note that the boundaries of these cycles impose a bilinear structure on the group $\calZ_0 (\Diamond_1)/\partial \calZ_1 (\Diamond_2)= ACH_0 (k,1)$ inducing an isomorphism $ACH_0 (k,1) \simeq \Omega_{k/\bbZ} ^1$ (see Theorem 6.4 in \cite{BE2}). We have a similar lemma that will be used in $\S 2$.
\begin{lemma}[\emph{c.f.} (6.22) in \cite{BE2}]\label{bilinearity2}Suppose that $k$ is algebraically closed. Consider the following set-theoretic map 
\begin{equation}f : k \times k^{\times} \to \calZ_0 (\Diamond_1)\end{equation}
$$(a,b) \mapsto \tuborg \left( \frac{1}{a}, b \right), & \mbox{ if } a \not = 0, \\ 0, & \mbox{ if } a = 0. \sluttuborg$$ Then, the map $f$ descends to a homomorphism
$$\overline{f}: k \otimes_{\bbZ} k^{\times}  \to \calZ_0 (\Diamond_1) / \left<\partial C_1 ^{\cdot},\partial C_2 ^{\cdot} \right>,$$ where $k$ is regarded as the additive abelian group and $k^{\times}$ is regarded as the multiplicative abelian group. This map $\overline{f}$ is in fact an isomorphism.
\end{lemma}

\begin{proof}That the set theoretic map $f$ descends to a homomorphism $\overline{f}$ follows immediately from the Lemma \ref{bilinearity}. That this gives an isomorphism can be seen as follows. As $k$ is algebraically closed, all generators of the free group $\calZ_0 (\Diamond_1)$ are $k$-rational. Define a homomorphism
$$g: \calZ_0 (\Diamond_1) \to k \otimes _{\bbZ} k^{\times}$$
$$\left(\frac{1}{a}, b\right) \mapsto a \otimes b.$$
By the Lemma \ref{bilinearity} again, the map $g$ descends to $$\overline{g}: \calZ_0 (\Diamond_1) / \left< \partial C_1 ^{\cdot}, \partial C_2 ^{\cdot} \right> \to k \otimes_{\bbZ} k^{\times}.$$ It is easy to see that $\overline{g}$ and $\overline{f}$ are inverse to each other.
\end{proof}

\subsection{Nontriviality of $ACH_1 (k,2)$}

Define two parametrized $1$-cycles in $\calZ_1 (\Diamond_2)$ $$\tuborg \Gamma_1 = \left\{ \left( t, t, \frac{ \left( 1- \frac{1}{2} t \right)^2}{1-t} \right) | t \in k \right\}, \\ 
\Gamma_2 = \left\{ \left( t, 1+ \frac{t}{6} , 1- \frac{t^2}{4} \right) | t \in k \right\}. \sluttuborg$$
The cycle $\Gamma_2$ is a variation of the cycle $\calZ_2$ that appeared in \S 6 of \cite{BE2}. The point is that modulo the boundaries of $C_1 ^{\cdot}$ and $C_2 ^{\cdot}$, the boundaries of $\Gamma_1$ and $\Gamma_2$ are equivalent (\emph{i.e.} $\partial \Gamma_1 \equiv \partial \Gamma_2 \mod \left< \partial C_1 ^{\cdot}, \partial C_2 ^{\cdot} \right>$) but $\Gamma_1$ and $\Gamma_2$ have distinct regulator values as we will see in the following two lemmas.

\begin{lemma}\label{lemma gamma one} \begin{enumerate}

\item $R_2 (\Gamma_1) = \frac{1}{4}$.
\item $ \partial \Gamma_1 \equiv (1,2) \mod \left< \partial C_1 ^{\cdot}, \partial C_2 ^{\cdot} \right>$.

More precisely, for the cycle $\overline{\Gamma}_1:= \Gamma_1 + C_1 ^{\left( \frac{1}{2}, \frac{1}{2} \right), 2}$, we have $R_2 \left( \overline{\Gamma}_1 \right) = \frac{1}{4}$ and $\partial \left( \overline{\Gamma}_1 \right) = (1,2)$.
\end{enumerate}
\end{lemma}

\begin{proof}(1) We have $\Gamma_1 \in \calM^2 (t_2)$ and $$-\nu^* \left( \frac{1-t_2}{x^3} \frac{dt_1}{t_1} \right) = - \frac{ \frac{1}{4} }{x (1-x)} \frac{dx}{x}  = \frac{1}{4} \cdot \frac{1}{x^2} \left( 1 + x + x^2 + \cdots \right) dx,$$ so that $R_2 (\Gamma_1) = \res_{x=0} \left( \frac{1}{4} \frac{1}{x^2} (1+ x+ x^2 + \cdots) dx \right) = \frac{1}{4} $.

(2) We can compute it directly:
$$\tuborg \partial_1 ^0 \Gamma_1 = (0, 1) = 0, \\
\partial _1 ^{\infty} \Gamma_1 = (\infty, \infty) = 0, \\
\partial_2 ^{0} \Gamma_1 = (2,2) + (2,2) = 2 \cdot (2,2) ,\\
\partial _2 ^{\infty} \Gamma_1 = (1,1) = 0. \sluttuborg$$
Since $\partial C_1 ^{\left( \frac{1}{2}, \frac{1}{2} \right), 2} = - 2 \cdot (2,2) + (1,2)$ and $R_2 (C_1 ^{\cdot}) = 0$, the assertions follow.
\end{proof}

\begin{lemma}\label{lemma gamma two}\begin{enumerate}
\item $R_2 (\Gamma_2) = - \frac{1}{24}$.
\item $\partial \Gamma_2 \equiv (1,2) \mod \left< \partial C_1 ^{\cdot}, \partial C_2 ^{\cdot} \right>$.

More precisely, if we let \begin{eqnarray*} \overline{\Gamma}_2 &:=& \Gamma_2 + \left\{ C_1 ^{\left( - \frac{1}{2}, \frac{1}{2} \right), \frac{2}{3} } + 3 C_1 ^{\left( - \frac{1}{3}, - \frac{1}{3} \right), -1 } - C_1 ^{\left( - \frac{1}{6}, - \frac{1}{6} \right), 2} \right. \\ &-& \left. C_1 ^{\left( - \frac{1}{6}, - \frac{1}{3} \right), 2} - C_1 ^{\left( - \frac{1}{2}, \frac{1}{2} \right) , 2 } + C_1 ^{\left( \frac{1}{2} , \frac{1}{2} \right), 2 } \right\} \\
& +& \left\{ - C_2 ^{ - \frac{1}{6}, \left( 4, -2 \right) } - C_2 ^{ - \frac{1}{6}, \left( -2, -2 \right) } + C_2 ^{\frac{1}{2} \left( \frac{2}{3} , \frac{3}{2} \right) } \right. \\ \ &-& \left. 3 C_2 ^{ - \frac{1}{6}, \left( 2, -1 \right)} - 3 C_2 ^{- \frac{1}{3}, \left( -1 , -1 \right) } + C_2 ^{\frac{1}{2} , \left( \frac{4}{3}, \frac{3}{2} \right) } \right\},\end{eqnarray*}
we have $R_2 \left( \overline{\Gamma}_2  \right) = - \frac{1}{24}$ and $\partial \overline{\Gamma}_2  = (1,2)$.
\end{enumerate}
\end{lemma}

\begin{proof}
(1) We have $\Gamma_2 \in \calM^2 (t_2)$ so that we use $-\nu^* \left( \frac{1-t_2}{x^3} \frac{dt_1}{t_1} \right)= - \left( \frac{ \frac{t^2}{4}}{t^3} \frac{dt}{t+6} \right)  = - \frac{dt}{4t (t+6)}$. Hence $R_2 (\Gamma_2) = \res_{t=0} \left( - \frac{dt}{4t (t+6)} \right) = - \frac{1}{24}$.

(2) By a direct computation, we have $$\tuborg \partial_1 ^0 \Gamma_2 = (-6, -8), \\
\partial _1 ^{\infty} \Gamma_2 = 0, \\
\partial _2 ^0 \Gamma_2 = \left( 2, \frac{4}{3} \right) + \left( -2, \frac{2}{3} \right), \\
\partial _2 ^{\infty} \Gamma_2 = 0, \sluttuborg$$ so that $\partial \Gamma _2 = - (-6, -8) + \left( 2, \frac{4}{3} \right) + \left( -2, \frac{2}{3} \right).$ Now, modulo $ \left< \partial C_1 ^{\cdot}, \partial C_2 ^{\cdot} \right>$ we prove that $\partial \Gamma_2 \equiv (1,2)$. Indeed,
$$ \tuborg - \partial C_2 ^{- \frac{1}{6}, \left( 4, -2 \right)} - (-6, -8) = - (-6, 4) - (-6, -2),\\
 - \partial C_2 ^{- \frac{1}{6}, (-2, -2)} - (-6, 4) = - (-6, -2) - (-6, -2),\\
\partial C_1 ^{\left( - \frac{1}{2}, \frac{1}{2} \right), \frac{2}{3} } + \left( -2, \frac{2}{3} \right) = - \left( 2, \frac{2}{3} \right),\\
\partial C_2 ^{\frac{1}{2} , \left( \frac{2}{3}, \frac{3}{2} \right)} - \left( 2, \frac{2}{3} \right) = \left( 2, \frac{3}{2} \right),\\
- 3 \partial C_2 ^{- \frac{1}{6}, (2, -1) } - 3 \left( -6, -2 \right) = -3 \left( -6 , 2 \right) - 3 \left( -6, -1 \right),\\
3 \partial C_1 ^{\left( - \frac{1}{3}, - \frac{1}{3} \right), -1} - 3 (-6, -1) = -6 (-3, -1), \\
-3 \partial C_2 ^{- \frac{1}{3} , \left( -1 , -1 \right) } - 6 (-3, -1) = 0,\\
 \partial C_2 ^{\frac{1}{2} , \left( \frac{4}{3}, \frac{3}{2} \right)} + \left( 2, \frac{4}{3} \right) + \left( 2, \frac{3}{2} \right) = \left( 2,2 \right),\\
 - \partial C_1 ^{\left( - \frac{1}{6}, - \frac{1}{6} \right), 2} -2 (-6, 2) = - (-3, 2),\\
- \partial C_1 ^{\left( - \frac{1}{6}, - \frac{1}{3} \right), 2} - (-6, 2) - (-3, 2) = - (-2, 2),\\
- \partial C_1 ^{\left( - \frac{1}{2}, \frac{1}{2}\right), 2} - (-2, 2) = (2,2),\\
\partial C_1 ^{\left( \frac{1}{2}, \frac{1}{2} \right), 2} + 2 (2,2) = (1,2).\sluttuborg $$
Since $R_2 (C_1 ^{\cdot}) = R_2 (C_2 ^{\cdot}) = 0$, we have $\partial \overline{\Gamma}_2 = (1,2)$ and $R_2 \left( \overline{\Gamma}_2 \right) = - \frac{1}{24}$.\end{proof}

\begin{prop}Let $\Gamma_3 := \overline{\Gamma}_1  - \overline{\Gamma}_2$. Specifically, 
\begin{eqnarray*} \Gamma_3 &=& \left\{ \left( x, x, \frac{ \left( 1- \frac{1}{2} x \right)^2}{1-x } \right) | x \in k  \right\} - \left\{ \left( x, 1 + \frac{x}{6} , 1- \frac{x^2}{4} \right) | x \in k \right\} \\
 &+& \left( - C_1 ^{\left( - \frac{1}{2}, \frac{1}{2} \right), \frac{2}{3} } -3 C_1 ^{\left( - \frac{1}{3} , - \frac{1}{3} \right), -1} + C_1 ^{ \left( - \frac{1}{6}, - \frac{1}{6} \right), 2} \right. \\ &+& \left. C_1 ^{\left( - \frac{1}{6}, - \frac{1}{3} \right), 2} + C_1 ^{\left( - \frac{1}{2}, \frac{1}{2} \right) , 2 } \right) \\ 
&+& \left( C_2 ^{- \frac{1}{6} , (4, -2)} + C_2 ^{- \frac{1}{6}, (-2, -2)} - C_2 ^{\frac{1}{2}, \left( \frac{2}{3}, \frac{3}{2} \right)}\right. \\ & +& \left. 3 C_2 ^{- \frac{1}{6}, (2, -1)} + 3 C_2 ^{- \frac{1}{3} , (-1, -1)} - C_2 ^{\frac{1}{2} , \left( \frac{4}{3}, \frac{3}{2} \right)} \right).\end{eqnarray*} Then, the cycle $\Gamma_3$ satisfies $\partial \Gamma_3 = 0$ and $R_2 \left( \Gamma_3 \right) = \frac{7}{24} \not = 0$.
\end{prop}

\begin{proof}It follows from Lemmas \ref{lemma gamma one} and \ref{lemma gamma two}: $\partial (\Gamma_3) = \partial \left( \overline{\Gamma}_1  - \overline{\Gamma}_2  \right) = (1,2) - (1,2) = 0$ and $R_2 (\Gamma_3) = R_2 \left( \overline{\Gamma}_1 \right) - R_2 \left( \overline{\Gamma}_2  \right) = \frac{1}{4} + \frac{1}{24} = \frac{7}{24} \not = 0$.
\end{proof}

As a corollary, we have:
\begin{thm}\label{theorem 1}When $k$ is a field of characteristic $0$, the above $1$-cycle $\Gamma_3$ represents a nontrivial class in $ACH_1 (k,2)$. In particular $ACH_1 (k,2) \not = 0$.
\end{thm}

It is conjecturally believed that this group is isomorphic to $k$. See \cite{P1} for a little more details about it.

\section{The additive dilogarithm}

Until the last, we let $k$ be an algebraically closed field of characteristic zero. As mentioned in the introduction, we suppose we have an isomorphism through the regulator map $R_2$:
\begin{equation}\label{star} R_2 : ACH_1 (k,2) \overset{\sim}{\to} k.\end{equation}The author doesn't have a proof of this assumption yet.

The $K$-theoretic version of \eqref{cathelineau sequence} discussed in \cite{BE2}
$$0 \to k \to TB_2 (k) \overset{\partial}{\to} k \otimes k^{\times} \to \Omega_{k/\bbZ} ^1 \to 0$$ had important classes of elements in $TB_2 (k)$ denoted by $\left< a \right>$ for $a \in k \backslash \{ 0, 1 \}$ with the properties
\begin{equation}\label{K generator}
\tuborg \rho\left( \left< a \right> \right) = a (1-a) \in k, \\ \partial \left( \left< a \right> \right) = 2 \left( a \otimes a + (1-a) \otimes (1-a) \right) \in k \otimes k^{\times},\sluttuborg
\end{equation}where $\rho: TB_2 (k) \to k$ is a homomorphism defined in Proposition 2.3 in \cite{BE2}. This group $TB_2 (k)$ is isomorphic to $k \oplus \beta(k)$ via $\rho \oplus \partial$ (see Lemma 3.7 in \cite{BE2}) as $k^{\times}$-modules, thus it is identified with $T\calP(k)$ of \eqref{cathelineau sequence}.

Our description using cycles begins with the additive Chow complex
$$\cdots \to \calZ_2 (\Diamond_3) \overset{\partial_2}{\to} \calZ_1 (\Diamond_2) \overset{\partial_1}{\to} \calZ_0 (\Diamond_1) \to 0.$$It induces the exact sequence
\begin{equation}\label{cycle sequence}0 \to \frac{\ker \partial_1}{{\rm im} \partial_2} \to \frac{\calZ_1 (\Diamond_2)}{ {\rm im} \partial_2 + \left< C_1 ^{\cdot}, C_2 ^{\cdot} \right>} \overset{\overline{\partial}_1}{\to} \frac{\calZ_0 (\Diamond_1)}{\left< \partial C_1 ^{\cdot}, \partial C_2 ^{\cdot} \right>} \to \frac{\calZ_0 (\Diamond_1)}{{\rm im} \partial_1} \to 0\end{equation} where $C_1 ^{\cdot}, C_2 ^{\cdot}$ are the $1$-cycles defined in \S 1. This sequence is equivalent to
$$0 \to k \to T\calP^{cy} (k) \overset{\overline{\partial}_1}{\to} k \otimes k^{\times} \to \Omega_{k/\bbZ} ^1 \to 0$$ by \eqref{star}, Lemma \ref{bilinearity2} and Theorem 6.4 in \cite{BE2}, where $T\calP ^{cy} (k)$ is the second group of the \eqref{cycle sequence}. Given the important roles of $\left< a \right>$ in $TB_2 (k)$, we may look for $1$-cycles $C_a \in \calZ_1 (\Diamond_2)$ with the analogous properties as \eqref{K generator}
\begin{equation}\label{cycle generator}
\tuborg R_2 \left( C_a \right) = a (1-a) \in k, \\ \overline{\partial}_1 \left( C_a \right) = \alpha * \left( \left( \frac{1}{a} , a \right) + \left( \frac{1}{1-a} , 1-a \right) \right)  \mbox{ for some } \alpha \in k^{\times}, \sluttuborg
\end{equation}where the last expression corresponds to $\alpha \cdot \left( a \otimes a + (1-a) \otimes (1-a) \right)$ in $k \otimes k^{\times}$ via the Lemma \ref{bilinearity2}. The definition of $C_a$ is given in \eqref{C_a}, but in any case its existence is more important so that we proceed to prove the results leaving its definition aside for a while. 

Recall from Proposition 6 in \cite{C} that $\beta (k)$ is the kernel of the map $k \otimes k^{\times} \to \Omega_{k/\bbZ} ^1$ mapping $a \otimes b \mapsto a db/b$.

\begin{lemma}\label{cycle identification}Under \eqref{star}, the map $R_2 \oplus \overline{\partial}_1 : T\calP ^{cy} (k) \to k \oplus \beta(k)$ is an isomorphism.
\end{lemma}
\begin{proof}Since the map $R_2 : ACH_1 (k) \to k$ gives an isomorphism, via this identification we have a splitting
$k \hookrightarrow T\calP^{cy} (k) \overset{R_2}{\to} k.$ The cokernel of $k \hookrightarrow T\calP^{cy} (k)$ is $\beta(k)$ by the exact sequence \eqref{cycle sequence}. This finishes the proof.\end{proof}

\begin{cor}Under \eqref{star}, we have identifications of the following three groups:
\begin{enumerate}
\item $T\calP(k)$ in \eqref{cathelineau sequence},
\item $TB_2 (k)$ in \cite{BE2},
\item $T\calP ^{cy} (k)$.
\end{enumerate}
\end{cor}

We remark now that the cycles $C_a$ generate the group $T\calP^{cy} (k)$ as a $k^{\times}$-module and they satisfy the Cathelineau identity.

\begin{lemma}Under \eqref{star}, every element in $T\calP^{cy} (k)$ can be written as a sum $\sum c_i * C_{a_i}$. In other words, $T\calP ^{cy}(k)$ is generated as $k^{\times}$-module by $C_a$'s.\end{lemma}

\begin{proof}The proof is essentially identical to that of Lemma 3.3 in \cite{BE2}. The properties \eqref{cycle generator} imply the desired result.
%
\end{proof}
\begin{lemma}Under \eqref{star}, the cycles $\{ C_a \}$ for $a \in k \backslash \{ 0, 1 \}$ satisfy the Cathelineau identity
$$C_a - C_b + a * C_{\frac{b}{a}} + (1-a)* C_{\frac{1-b}{1-a}} \equiv 0 \ \ \mbox{ in } T\calP^{cy} (k).$$
\end{lemma}
\begin{proof}Let $D_{a,b}$ be the left hand side of the expression. Since $R_2 (C_a) = a (1-a)$, we have $$R_2 \left(C_a - C_b + a* C_{\frac{b}{a}} + (1-a) * C_{\frac{1-b}{1-a}}\right)$$
$$= a (1-a) - b (1-b) + a^3 \cdot \frac{b}{a} \left( 1 - \frac{b}{a} \right) + (1-a)^3 \cdot  \left( \frac{1-b}{1-a} \right) \left( 1 - \frac{1-b}{1-a} \right) = 0.$$ Since $\overline{\partial}_1 (C_a) = \alpha \cdot \left( a \otimes a + (1-a) \otimes (1-a) \right) $ in $k \otimes k^{\times}$ through the identification of Lemma \ref{bilinearity2}, by Proposition 6 of \cite{C} we have $\overline{\partial}_1 (D_{a,b}) = 0$ in $\beta(k) \subset k \otimes k^{\times}$. Hence by the Lemma \ref{cycle identification}, the cycle $D_{a,b}$ must represent the zero class in $T \calP^{cy} (k)$.\end{proof}

The summary of the above discussion is the following theorem:
\begin{thm}\label{dilog}Let $k$ be an algebraically closed field of characteristic $0$. Assume \eqref{star} that the regulator $R_2: ACH_1 (k,2) \to k$ gives an isomorphism. Then, we have the additive $4$-term motivic exact sequence
$$0 \to ACH_1 (k,2) \to T\calP^{cy} (k) \to k \otimes k^{\times} \to ACH_0 (k,1) \to 0$$ obtained from the additive higher Chow complex. In addition, we have $ACH_1(k,2) \simeq k$, $ACH_0 (k,1) \simeq \Omega_{k/\bbZ} ^1$, $T\calP^{cy}(k) \simeq k \oplus \beta(k)$ and
$$\beta (k) := \frac{k [ k \backslash \{ 0, 1 \} ]}{\left(\left< a \right> - \left< b \right> + a \left< \frac{b}{a} \right> + (1-a) \left< \frac{1-b}{1-a} \right>\right) }.$$ There are classes $\{C _a \}$ in $T\calP^{cy} (k)$ for $a \in k \backslash \{ 0, 1 \}$ represented by $1$-cycles in $\calZ_1 (\Diamond_2)$ that generate $T\calP^{cy} (k)$ as a $k^{\times}$-module and satisfy the Cathelineau identity
$$C_a - C_b + a * C_{\frac{b}{a}} + (1-a)* C_{\frac{1-b}{1-a}} \equiv 0 \ \ \mbox{ in } T\calP^{cy} (k).$$ \end{thm}

\bigskip The rest of the section is devoted in writing down the cycles $C_a$ satisfying \eqref{cycle generator} concretely. The definition is given in \eqref{C_a}. They are variations of the cycle $Z(1-2a)$ in the last section of \cite{BE2}. We modify this cycle to equip a better property.

Let $Q (a) := \left\{ \left( t, 1 + \frac{t}{2}, 1- \frac{a^2 t^2}{4} \right) | t \in k \right\}$.

\begin{lemma}\begin{enumerate}
\item $R_2 \left( Q (1-2a) \right) =  - \frac{1}{2} \left( a (1-a) \right) - \frac{1}{8}$.
\item $\partial Q \left( 1-2a \right) \equiv \left( \frac{1}{a}, a \right) + \left( \frac{1}{1-a}, 1-a \right) + (1,2) \mod \left< \partial C_1 ^{\cdot}, \partial C_2 ^{\cdot} \right>.$

More precisely, for the cycle

\begin{eqnarray*}\widetilde{Q}(a)&= &Q(1-2a) + \left\{ C_1 ^{\left( \frac{1-2a}{2} , - \frac{1-2a}{2} \right) 1- \frac{1}{1-2a}} + C_1 ^{\left( - \frac{1}{4}, - \frac{1}{4} \right), -1 } \right. \\ 
&+& \left. C_1 ^{\left( \frac{1}{2}, \frac{1-2a}{2} \right), 2-2a} - C_1 ^{\left( - \frac{1}{2}, \frac{1}{2} \right), 2-2a}  - C_1 ^{\left( - \frac{1}{2}, \frac{1-2a}{2} \right), -2a} \right. \\ 
&+ & \left. C_1 ^{\left( a, 1-a \right), 2} - C_1 ^{\left( -a, a \right), -2a} - C_1 ^{\left( \frac{a}{2}, \frac{a}{2} \right), -1 } \right\} \\ 
& +& \left\{ C_2 ^{\frac{1-2a}{2} , \left( 1- \frac{1}{1-2a}, \frac{1-2a}{-2a} \right)} + C_2 ^{\frac{1-2a}{2}, \left( 1+ \frac{1}{1-2a}, \frac{1-2a}{-2a} \right)} \right. \\
& +& \left. C_2 ^{\frac{1-2a}{2}, \left( 2-2a, \frac{1}{-2a} \right)} +   C_2 ^{- \frac{1}{2}, \left( -1, -2a \right)} + C_2 ^{- \frac{1}{4} , (-1, -1)} \right. \end{eqnarray*} \begin{eqnarray*}
&-&  \left. C_2 ^{- \frac{1}{2} , \left( 2a, 2-2a \right) }   -  C_2 ^{\frac{1-2a}{2} , \left( \frac{1}{-2a}, -2a \right)} + C_2 ^{a, \left( a, -2 \right)} \right. \\
&+& \left.  C_2 ^{1-a, \left( 1-a, 2 \right)} +   C_2 ^{a, (2, -1 )} - C_2 ^{\frac{a}{2}, \left( -1, -1 \right)} \right\},\end{eqnarray*}
we have $R_2 (\widetilde{Q} (a)) = - \frac{1}{2} \left( a (1-a) \right) - \frac{1}{8}$ and $$\partial \left( \widetilde{Q}(a)\right) = \left( \frac{1}{a} , a \right) + \left( \frac{1}{1-a}, 1-a \right) + (1,2).$$
\end{enumerate}
\end{lemma}

\begin{proof}(1) Since $Q(a) \in \calM^2 (t_2)$ we use $- \nu^* \left( \frac{1-t_2}{x^3} \frac{dt_1}{t_1} \right) = - \frac{a^2}{4t} \frac{dt}{t+2}$ so that $$R_2 \left( Q (a) \right) = \res_{t=0} \left( - \frac{a^2}{4t} \frac{dt}{t+2} \right) = - \frac{a^2}{8}.$$ Thus, the value $R_2 \left( Q(1-2a) \right) = - \frac{ (1-2a)^2}{8} = - \frac{1}{2} \left( a (1-a) \right) - \frac{1}{8}$.

(2) Notice that $$\tuborg \partial_1 ^0 Q (a) = (-2, 1-a^2), \\
\partial_1 ^{\infty} Q (a) = 0, \\
\partial _2 ^0 Q (a) = \left( \frac{2}{a} , 1 + \frac{1}{a} \right) + \left( - \frac{2}{a}, 1- \frac{1}{a} \right), \\
\partial_2 ^{\infty} Q (a) = 0, \sluttuborg $$
so, the cycle $\partial Q (a) = - \left( -2, 1-a^2 \right) + \left( \frac{2}{a} , 1 + \frac{1}{a} \right) + \left( - \frac{2}{a} , 1 - \frac{1}{a} \right).$ 

Now, we have
\begin{eqnarray*}  - \partial C_2 ^{ - \frac{1}{2}, (1-a, 1+a) } - (-2, 1-a^2 ) &=& - (-2, 1-a) - (-2, 1+a), \\
 \partial C_1 ^{\left( \frac{a}{2} , - \frac{a}{2} \right) , 1- \frac{1}{a} } + \left( - \frac{2}{a}, 1- \frac{1}{a} \right) &=& - \left( \frac{2}{a}, 1 - \frac{1}{a} \right),\\
 \partial C_2 ^{\frac{a}{2} , \left( 1- \frac{1}{2}, \frac{a}{a-1} \right)} - \left( \frac{2}{a} , 1- \frac{1}{a} \right) &=& \left( \frac{2}{a}, \frac{a}{a-1} \right), \end{eqnarray*} \begin{eqnarray*}
 \partial C_2 ^{\frac{a}{2} , \left( 1+ \frac{1}{a} , \frac{a}{a-1} \right)} + \left( \frac{2}{a} , 1 + \frac{1}{a} \right) + \left( \frac{2}{a} , \frac{a}{a-1} \right) &=& \left( \frac{2}{a}, \frac{a+1}{a-1} \right), \\
 \partial C_2 ^{\frac{a}{2} , \left( a+1, \frac{1}{a-1} \right)} + \left( \frac{2}{a} , \frac{a+1}{a-1} \right)& =& \left( \frac{2}{a} , a+1 \right) + \left( \frac{2}{a} , \frac{1}{a-1} \right),\\
 - \partial C_2 ^{\frac{a}{2} , \left( \frac{1}{a-1} , a-1 \right) } + \left( \frac{2}{a} , \frac{1}{a-1} \right) &=& - \left( \frac{2}{a} , a-1 \right),\\
 \partial C_1 ^{\left( - \frac{1}{2} , \frac{1}{2} \right), 1+a } - (-2, 1+a) &=& (2, 1+a),\\
 \partial C_2 ^{- \frac{1}{2} , \left(-1, a-1 \right)} - (-2, 1-a)& =& - \left( -2, -1 \right) - (-2, a-1),\\
\partial C_1 ^{\left( - \frac{1}{4}, - \frac{1}{4} \right) , -1 } - (-2, -1)& =& -2 (-4, -1),\\
\partial C_2 ^{- \frac{1}{4} , (-1, -1)} - 2 (-4, -1) &=& 0,\\
- \partial C_1 ^{\left( - \frac{1}{2} , \frac{a}{2} \right), a-1 } - (-2, a-1) - \left( \frac{2}{a} , a-1 \right) &=& - \left( \frac{2}{a-1} , a-1 \right),\\
\partial C_1 ^{\left( \frac{1}{2} , \frac{a}{2} \right), a+1 } + (2, 1+a) + \left( \frac{2}{a} , 1+a \right) &=& \left( \frac{2}{a+1} , a+1 \right). \end{eqnarray*}
Let $Q' (a)$ be the sum of all $C_1 ^{\cdot}$'s and $C_2 ^{\cdot}$'s of the above equations. Addition of all of the above equations give lots of cancellations and we end up with
$$ \partial \left( Q (a) + Q' (a) \right)= - \left( \frac{2}{a-1}, a-1 \right) + \left( \frac{2}{a+1}, a+1 \right).$$ Plugging in $1-2a$ in the place of $a$, we have
$$\partial \left( Q (1-2a) + Q' (1-2a) \right) = - \left( - \frac{1}{a}, - 2a \right) + \left( \frac{1}{1-a} , 2 (1-a) \right),$$ and from equations
\begin{eqnarray*} \tuborg - \partial C_1 ^{ (-a, a), -2a} - \left( - \frac{1}{a}, -2a \right) = \left( \frac{1}{a}, -2a \right),\\
\partial C_2 ^{a, (a, -2)} + \left( \frac{1}{a}, - 2a \right) = \left( \frac{1}{a}, a \right) + \left( \frac{1}{a} , -2 \right),\\
 \partial C_2 ^{1-a, \left( 1-a, 2 \right)} + \left( \frac{1}{1-a} , 2 (1-a) \right) = \left( \frac{1}{1-a}, 1-a \right) + \left( \frac{1}{1-a} , 2 \right),\\
\partial C_2 ^{a, (2, -1)} + \left( \frac{1}{a} , -2 \right) = \left( \frac{1}{a}, 2 \right) + \left( \frac{1}{a}, -1 \right),\\
\partial C_1 ^{\left( \frac{a}{2}, \frac{a}{2} \right), -1 } + \left( \frac{1}{a}, -1 \right) = 2 \left( \frac{2}{a}, -1 \right),\\
- \partial C_2 ^{\frac{a}{2} , \left( -1, -1 \right)} + 2 \left( \frac{2}{a}, -1 \right) = 0,\\
\partial C_1 ^{(a, 1-a) ,2 } + \left( \frac{1}{a}, 2 \right) + \left( \frac{1}{1-a}, 2 \right) = (1,2), \sluttuborg \end{eqnarray*}
we obtain $\partial \widetilde{Q} (a) = \left( \frac{1}{a}, a \right) + \left( \frac{1}{1-a} , 1-a \right) + (1,2)$.
\end{proof}

From the above, we see that $\partial \left( \widetilde{Q}(a) - \overline{\Gamma}_1 \right) = \left( \frac{1}{a} , a \right) + \left( \frac{1}{1-a} , 1-a \right)$. But, it does not have a right regulator value: $R_2 \left(\widetilde{Q}(a)  - \overline{\Gamma}_1 \right) = - \frac{1}{2} \left( a (1-a) \right) - \frac{1}{8} - \frac{1}{4} = - \frac{1}{2} \left( a (1-a) \right) - \frac{3}{8}$. We remedy this situation.

Let $\alpha= \sqrt[3]{-2}$ and $\alpha' = \sqrt[3]{\frac{-18}{7}}$. Then, $R_2 \left( \alpha * \left( \widetilde{Q}(a) - \overline{\Gamma}_1  \right) \right) = a (1-a) - \frac{3}{4}$. Since $R_2 \left( \Gamma_3 \right) = \frac{7}{24}$ and $\partial \Gamma_3 = 0$, if we let \begin{equation}\label{C_a}C_a := \alpha * \left( \widetilde{Q} (a) - \overline{\Gamma}_1  \right) - \alpha' * \Gamma_3, \end{equation} then as we desired we have
\begin{equation}\label{Cathelineau element}\tuborg R_2 (C_a) = a (1-a) - \frac{3}{4} + \frac{3}{4} = a (1-a), \\
\partial C_a = \partial \left( \alpha * \left( \widetilde{Q}(a)- \overline{\Gamma}_1 \right) \right) = \alpha * \left( \left( \frac{1}{a}, a \right) + \left( \frac{1}{1-a} , 1-a \right) \right). \sluttuborg \end{equation}

\begin{remark2}{\rm In connection with the third problem of D. Hilbert (see \cite{D1}), several authors observed interesting similarities between the additive motivic exact sequence \eqref{cathelineau sequence} and the basic exact sequence for the scissors congruence group of the $3$-dimensional Euclidean space (see \cite{BE2}, \cite{G3}). A general discussion on scissors congruence can be found in \cite{D3}. In this analogy, one may regard the regulator map $R_2$ as the \emph{volume map}, and the boundary map $\overline{\partial}_1$ as the \emph{Dehn invatiant} map. Certainly the regulator $R_2$ satisfies the property $R_2 (\alpha * C) = \alpha^3 R_2 (C)$ as seen in Remark \ref{volume}, and one may wish to interpret this $*$-action of $k^{\times}$ as the enlargement by $\times \alpha$ in the $3$-dimensional space. As observed by Sydler in \cite{S}, a class in the scissors congruence group is determined by its volume and the Dehn invariant, just like our group $T\calP^{cy}(k)$ is determined by the images of $R_2$ and $\overline{\partial}_1$ as seen in Lemma \ref{cycle identification}. However it is still mysterious to the author why this interesting phenomena occur, and how one can associate some polyhedra to cycles.

}
\end{remark2}

\vskip0.5cm

\noindent \textbf{Acknowledgment} This paper is based on a chapter of the author's doctoral thesis at the University of Chicago. The author would like to thank Spencer Bloch, Jean-Louis Cathelineau, H\'el\`ene Esnault, Madhav Nori, Kay R\"ulling and the referee.

\end{document}